\documentclass{article}
\usepackage{amssymb, amsmath, amsfonts}
\begin{document}

\newtheorem{theo}{Theorem}
\newtheorem{lem}{Lemma}

\begin{center}
{\huge {\bf Generalization of the Critical Volume NTCP\\
\vskip0,3cm Model in the Radiobiology}}
\end{center}
\vskip0,4cm

\begin{center}
{\Large {\bf Alexander Bulinski$^{a,}
\footnote{
This work is partially supported by the RFBR grant 03-01-00724, by the grant 1758.2003.1 of Scientific Schools
and by the INTAS grant 03-51-5018.
}
$, Andrei Khrennikov$^{b}$ }}
\end{center}
\vskip0,1cm

$^a$ {\it Dept. of Mathematics and Mechanics of the Moscow State University, Moscow 119992, Russia}\\
E-mail: bulinski@mech.math.msu.su

$^b$ {\it Int. Center for Mathematical Modelling in Physics and Cognitive Sciences, MSI, University of
V\"axj\"o, S-35195, Sweden}\\ E-mail: Andrei.Khrennikov@msi.vxu.se

\vskip0,5cm
{\large {\bf Abstract}}
\vskip0,3cm

A generalization of the well known critical volume NTCP model is proposed to take into account
dependence of the functional subunits of irradiated organ (or tissue).
A new statistical version of the CLT is established to analyze the corresponding random fields.

\vskip0,5cm
{\bf AMS} classification: 60F05, 62E20, 62G15, 62P10.

{\it Key words}: NTCP, dependence conditions, random fields, CLT,
statistical version of the CLT. \vskip1cm

{\large {\bf 1. Introduction}}
\vskip0,3cm

The problem of finding optimal radiation doses for organs or tissues in therapy of cancer belongs
to the principal ones in the modern Radiobiology (see, e.g., \cite{Smith}, \cite{JonDale},
\cite{LMHK}, \cite{MBD}, \cite{SSWF}, \cite{LDSS}, \cite{PPVSK} and references therein). The
complexity of this problem is related to nondeterministic character of the oncological therapy
results.

The aim of this paper is to study stochastic models for collective effects in the behaviour of the
irradiated cells.  We provide a generalization of the well-known {\it critical volume} (CV) {\it
normal tissue complication probability} (NTCP) model to comprise a concept of {\it functional
subunits} (FSUs) dependence. Such a model is beyond the scope of \cite{B} and for its investigation
new limit theorems are required. Note that here the links between Probability and Geometry are
stipulated by the dependence structure of a random field under consideration which is governed by
the configuration of a graph used as a parameter set. Moreover, it seems natural from the
biological view-point to assume some dependence in the collective performance of cells (and FSUs).

The description of the FSUs behaviour by means of non-binary random variables is considered as
well. Other important biological response models and further research directions are tackled in the
last Section.

\vskip1cm

{\large {\bf 2. Accuracy of the CV NTCP model}}
\vskip0,3cm

We recall the basic {\it critical volume} model (see \cite{Jack}, \cite{Nie}) and after that
consider more carefully its framework. The organ (or tissue) modelled is assumed to be composed of
independent FSUs and it is supposed that complications in its functioning arise only if
sufficiently many FSUs ("the functional reserve") are destroyed. More precisely (see, e.g.,
\cite{SNSG}), the assessment of the impact of irradiation is divided into two stages. The first one
is the reaction of the cells forming an FSU which gives the probability of death or survival for
the FSU. The second one is a probability to have sufficient number of FSUs survived the irradiation
in order to maintain the organ's functionality (it means the complications absence probability).
Note that the {\it serial} and {\it tumor control} models are special cases when the functional
reserve consists of one or all the FSUs respectively. One speaks also of the {\it critical element}
models. The last special models are natural for organs having the one-dimensional shape (e.g.,
spinal chord). For other type of organs it is desirable to consider the "parallel architecture".

For the sake of simplicity we start with a model of organ consisting of $n$ similar FSUs. Introduce
independent identically distributed (i.i.d.) random variables (r.v's) $X_1,\ldots,X_n$ defined on
some probability space $(\Omega, \mathcal{F},{\sf P})$ such that
\begin{equation}\label{1}
{\sf P}(X_i=1) =p,\;\;\;{\sf P}(X_i=0) = q\;\;{\mbox where}
\;\;q=1-p,\;\; 0<p<1.  \end{equation} More exactly consider an array
of r.v's $X_i = X_i(D)$ where $D$ is a positive parameter ({\it the
irradiation dose}). Assume that $X_i(D)$ shows the state of the
$i$-th FSU after irradiation of dose $D$ ($i=1,\ldots,n)$. Namely,
the random event $\{X_i =1\}$ means that the $i$-th FSU is killed
and $\{X_i=0\}$ corresponds to the case that the $i$-th FSU
survives. As usual the argument $\omega \in \Omega$ is omitted and
we write $X_i(D)$ (or simply $X_i$) instead of $X_i(\omega,D)$.

Set $S_n(D) = \sum_{i=1}^n X_i(D)$. In other words consider a random variable equal to the number
of FSUs killed due to irradiation dose $D$. Thus NTCP is ${\sf P}(S_n(D) \geq L)$ with {\it the
threshold} $L$ being some positive integer.

Using the convergence rate estimate in the central limit theorem (see, e.g.,
\cite{CT}, p. 323) one has
\begin{equation}\label{2} \sup_{-\infty<x<\infty} |{\sf
P}(S_n(D) \geq x) - 1 + \Phi(z)| \leq c/\sqrt{np(D)q(D)}
\end{equation}
where
$z = (x-np(D))/\sqrt{np(D)q(D)}$, $p(D)$ and $q(D)$ appear in \eqref{1} for
$X_i = X_i(D)$, positive constant $c \leq 0.7975$ and the c.d.f. of a standard normal
law
$$ \Phi(z) = \frac{1}{\sqrt{2\pi}}\int_{-\infty}^z e^{-\frac{u^2}{2}}
du,\;\;\;z \in \mathbb{R}.  $$

Thus if $np(D)q(D)$ is large enough, \footnote{Otherwise one has to use different approximations.}
for any threshold $x$ (possibly depending on $n$ and $p(D)$) the following
relation holds
\begin{equation}\label{3} {\sf P}(S_n(D) \geq x) \approx 1 -
\Phi(z)
\end{equation}
where $z$ is introduced above and the proximity is
evaluated by the right hand side of \eqref{2}. It is well known
that the power $-\frac{1}{2}$ of $n$ in \eqref{2} is the best possible
(see, e.g., \cite{Pet} for that and also about non uniform estimates of
the convergence rates in CLT for independent summands).

Consequently, for a given value $\gamma \in (0,1)$ (close enough to 1), using the table of function
$\Phi$
one can find $z = z_{\gamma}$ as the unique root of the equation $\Phi(z) =
\gamma$. Then the probability of complications is approximately $1-\gamma$ (with exactly specified
boundaries) and the threshold
\begin{equation}\label{4}
x_{\gamma} = np(D) + (np(D)q(D))^{1/2}z_{\gamma}.
\end{equation}
Evidently, $x_{\gamma} = x_{\gamma}(n,p(D))$.

The distribution of r.v. $S_n$ has atoms at the points $0,\ldots,n$.
Therefore one could choose an integer threshold
$L_{\gamma}=[x_{\gamma}]$ where $[\cdot]$ stands for the integer part of a
real number.  In this case a trivial estimate for the continuity module of
the function $\Phi$ shows that $$ |{\sf P}(S_n(D) \geq L_{\gamma}) - 1 +
\Phi(z)| \leq \frac{1}{\sqrt{np(D)q(D)}}\left(c+\frac{1}{\sqrt{2\pi}}\right).
$$
We
observe that if $np(D)q(D)$ is large enough then
a search for an integer threshold $L_{\gamma}$ is not important.

Note in passing that for the number of "successes" $S_n$ in the
Bernoulli scheme
(of $n$ independent trials with probability $p$ for "success") one can apply
(see, e.g., \cite{Weiss}) a little bit different approximation using for $0
\leq k \leq m \leq n$ the relation $$ {\sf P}(k \leq S_n \leq m) = \Phi(t_2) -
\Phi(t_1) + \frac{q-p}{6\sqrt{2\pi}\sigma} \{(1-t^2)e^{-t^2/2}\}|_{t_1}^{t_2} +
\Delta $$ where $\sigma = \sqrt{npq}$, $t_1 = (k-\frac{1}{2} - np)/\sigma$,
$t_2 = (m + \frac{1}{2} - np)/\sigma$ and the error term $\Delta$ satisfies for
$\sigma \geq 5$ the inequality $$ |\Delta| \leq (0.12 + 0.18|p-q|)\sigma^{-2} +
e^{-3\sigma/2}.  $$

Actually we deal with an equivalent description of the well-known critical volume model. Namely,
suppose that the volume \footnote{Possibly a length or an area, the interpretation depends on the
model of an organ} of the irradiated organ is $V$ and let $V_i$ represent the volume of the $i$-th
FSU ($i=1,\ldots,n)$.  Clearly instead of $S_n(D)$ (the number of killed FSUs) we could consider
the random damage volume $\widetilde{V} = \sum_{i: X_i = 1} V_i$. In the case when $V_i = V/n$,
$i=1,\ldots,n$, one has
\begin{equation}\label{5} \widetilde{V} = V\frac{S_n(D)}{n}.  \end{equation}
Thus we come to description of the irradiation result in terms of the damage
volume and one can specify the threshold $v_c$ for ${\sf P}(\widetilde{V} \geq
v_c)$.

Formula \eqref{5} suggests that it is natural to introduce a threshold of the type \footnote{We do
not use in this paper the theory of large deviations for sums of r.v's.} $x = \varkappa n$ where
$\varkappa \in (0,1)$ is the fraction of killed FSUs. Thus \begin{equation}\label{6} \varkappa = p
+ c(p(1-p))^{1/2}\;\;\;{\mbox where} \;\; p = p(D),\;\; c = z_{\gamma}n^{-1/2}.
\end{equation}

Note that  $z_{\gamma} \geq 0$ for $\gamma \geq 1/2$ and consequently $c \geq 0$. Evidently for
$c=0$ (i.e. $\gamma = 1/2$) one has $\varkappa = p$. For each $c>0$ the
graph of a function
$\varkappa = \varkappa(p)$
has the following features.
One can easily verify that $\varkappa(p_1) =1$ for $p_1 = 1/(1+c^2)$, $\varkappa'(0+) = +\infty$,
$\varkappa'(1-) = -\infty$ and the concave function $\varkappa(p)$ attains its maximum
$\varkappa_* = \frac{1}{2}(1+\sqrt{1+c^2})$ at the point $p_* = \frac{1}{2}(1
+ \frac{1}{\sqrt{1+c^2}})$. Moreover, for $\gamma \geq 1/2$ one has $$ 0 \leq
\inf_{0\leq p\leq 1} (\varkappa(p) - p) \leq \sup_{0\leq p \leq
1}(\varkappa(p)-p) \leq z_{\gamma}/(2\sqrt{n}).
$$
Note also that if $p(D) \geq p_1$ then relation \eqref{6} can not be
satisfied for any $\varkappa \in (0,1)$.

On the other hand, given $n \in \mathbb{N}$, $\gamma \geq 1/2$ (i.e. $c =
z_{\gamma}/\sqrt{n}$) and $\varkappa \in (0,1)$, there is a unique root $p =
\overline{p}$ of equation \eqref{6} \begin{equation}\label{7} \overline{p}=
\left(\varkappa + \frac{c^2}{2} + c\left(\varkappa - \varkappa^2 +
\frac{c^2}{4}\right)^{1/2}\right)/(1+c^2).
\end{equation}

We remark also that $p(D)$ should be
nondecreasing function on $(0,\infty)$ and if $p(D)$ is continuous then for any $\varkappa \in
(\varkappa_1,\varkappa_2)$ where $\varkappa_1 = \varkappa(\inf_{D>0} p(D))$ and $\varkappa_2 =
\min\{1,\varkappa(\sup_{D>0}p(D))\}$ there exists (unique if $p$ is
strictly increasing) $\overline{D}$ such that $p(\overline{D}) = \overline{p}$.

Now we discuss the models providing $p(D)$.
Assume that every FSU consists of $n_0$ cells. The {\it surviving fraction} of these cells after
irradiation of dose $D$ is determined (see, e.g., \cite{LDSS}) by
$$
SF(D) = \exp\{-\alpha D\}
$$
where $\alpha >0$ is the radiosensivity \footnote{One writes also $SF(D) = \exp\{-D/D_0\}$ where
$D_0$ is called the {\it mean lethal dose.}} of the cells. This is a
so-called
{\it single-hit model}. Suppose that each cell of FSU behaves in the same
manner as other ones. Usually one admits that an FSU can regenerate from a single
surviving cell, which means it is disabled only when no cell survives. Thus the
probability of killing an FSU  due to irradiation of dose $D$ is
\begin{equation}\label{8}
p(D) = (1 - e^{-\alpha D})^{n_0}.
\end{equation}
To obtain \eqref{8} one supposes that all $n_0$ cells in a FSU
evolve independently of each other.

Now assume that every cell contains $m$ targets, each of them must be hit
at least once to inactivate the cell. Then the probability that all targets of
a cell will be hit at least once is $(1 - e^{-\alpha D})^m$. Thus for this {\it
multi-target model} $$ SF(D) = 1 - (1-e^{-\alpha D})^m,\;\;\;p(D) = (1 -
e^{-\alpha D})^{mn_0}.  $$ Most experimental survival curves have an initial
slope whereas the multi-target/single-hit model predicts no initial slope. To
have a more adequate description one uses the family of functions $$ SF(D) =
e^{-\alpha D}(1 - (1 - e^{-\beta D})^m),\;\;\;\alpha >0,\;\;\beta>0.  $$ Note
(see, e.g., \cite{JonDale}) that for $SF(D)$ a {\it linear quadratic} (LQ)
model is also widely used with $$ SF(D) = e^{-(\alpha D + \beta
D^2)},\;\;\;\alpha>0,\;\;\beta>0.  $$

Now we concentrate on the generalization of the CV model considered above.

\vskip1cm

{\large {\bf 3. Variant of the central limit theorem for dependent random fields}}
\vskip0,3cm

Let $X(D) = \{X_j(D),j \in \mathbb{Z}^d\}$ $(d \geq 1)$ be a family of random fields defined on a
probability space $(\Omega,\mathcal{F},{\sf P})$ for $D>0$. Employing instead of an integer lattice
$\mathbb{Z}^d$ a parameter set $T = \delta \mathbb{Z}^d$, with $\delta >0$, one can easily
reformulate all the results for a family of FSUs assuming, e.g., that the $k$-th FSU is a cube with
a center at a point $t_k \in T$ and with an edge length equal to $\delta$.  In other words one can
use a scale appropriate to the problem under consideration. Thus without loss of generality we
restrict ourselves to the study of a random field $X(D)$ on a lattice $\mathbb{Z}^d$.  Moreover, we
can assume that the random variable $X_j(D)$ describes the state of the corresponding FSU after its
irradiation of dose $D$. This gives us a possibility to consider not only the death and survival of
an FSU but also to consider the "intermediate" states. Then the general (collective) effect of
irradiation is represented by the following sum
$$ S(U,D) = \sum_{j \in U} X_j(D) $$ where U is a finite subset of
$\mathbb{Z}^d$. For a fixed $D$ we also write simply $X_j$ and $S(U)$.

For a finite set $I \subset \mathbb{Z}^d$ with cardinality $|I|$ introduce the $\sigma$-algebra
$\mathcal{A}(I)= \sigma\{X_j, j \in I\}$, that is consider the $\sigma$-field generated by a field $X$
over a set $I$.

There are different methods (see, e.g., \cite{BuMSU}, \cite{Douk}) to
describe the dependence structure of a field $X$. Here we use the maximum
correlation coefficient for $\mathcal{A}(I)$ and $\mathcal{A}(J)$ over finite
disjoint sets $I,J \subset \mathbb{Z}^d$ which is defined as follows
\begin{equation}\label{9}
\rho(I,J) = \sup\{|corr(\xi,\eta)|: \;\xi \in L^2(\Omega,\mathcal{A}(I),{\sf P}),\;\; \eta \in
L^2(\Omega,\mathcal{A}(J),{\sf P})\}
\end{equation}
where $corr(\xi,\eta)$ is the correlation coefficient for (nondegenerate, square integrable)
real-valued random
variables $\xi$ and $\eta$ measurable with respect to $\sigma$-algebras
$\mathcal{A}(I)$ and $\mathcal{A}(J)$.

Assume that for all finite disjoint sets $I,J \subset \mathbb{Z}^d$, some positive $c_0$ and
$\lambda$ one has
\begin{equation}\label{10}
\rho(I,J) \leq c_0 |I||J| (dist(I,J))^{-\lambda}
\end{equation}
where
$$
dist(I,J) = \min\{\|q - j\|:\;q \in I,\;j\in J\},\;\;\; \|z\|= \max_{1\leq k
\leq d} |z_k|,\;z \in \mathbb{Z}^d.  $$

{\bf Remark 1.} Employing condition \eqref{10} has the following
motivation. In many stochastic models it is reasonable to assume
that dependence between the random variables $\{X_j, j \in I\}$ and
$\{X_j,j\in J\}$ is rather small if the distance between $I$ and $J$
is large enough. However, due to the paper \cite{Dob} it was
realized that for random fields (in contrast to stochastic processes
corresponding to the case $d=1$) one cannot, in general, measure the
dependence between $\mathcal{A}(I)$ and $\mathcal{A}(J)$ only in
terms of the distance between $I$ and $J$. Namely in many situations
the dependence between $\mathcal{A}(I)$ and $\mathcal{A}(J)$ could
increase for sets $I$ and $J$ growing, e.g., in such a way that the
distance between them is fixed. Dependence notions based on
correlations are quite familiar in various domains of applied
probability. Appearance of the factors $|U|$, $|V|$ and parameters
$c, \lambda$ in \eqref{10} is intended to account, in a qualitative
sense, for the effect of possible increase of dependence between
$\sigma$-algebras $\mathcal{A}(I)$ and $\mathcal{A}(J)$ when $I$ and
$J$ are growing so that the $dist(I,J)$ is preserved. Besides,
\eqref{10} implies the same (i.e. power-type) decrease of
correlations when $I$ and $J$ are moved apart so that the distance
between them tends to infinity. Of course, a simple and natural
hypothesis of $m$-dependence is a particular case of our condition.
Recall that a random field $X = \{X_j, j \in \mathbb{Z}^d\}$ is
$m$-dependent (with some $m >0)$ whenever $\mathcal{A}(I)$ and
$\mathcal{A}(J)$ are independent if $dist(I,J) \geq m$. Thus, we
include, in particular, a useful model of dependent nearest
neighbouring FSUs.  See also Remark 2 and Section 5.

First of all we establish the central limit theorem (CLT) with convergence rate
for partial sums \begin{equation}\label{11} S(U_n) = \sum_{j \in U_n} X_j,\;\;n
\in \mathbb{N}, \end{equation} of multi-indexed dependent r.v's where
summation is carried over the {\it integer cubes} $U_n = [-n,n]^d \cap
\mathbb{Z}^d$, $n \in \mathbb{N}$.

\begin{theo}\label{1} Let $X(D) = \{X_j(D), j \in \mathbb{Z}^d\}$, $D>0$,  be a family of
strictly stationary random fields
such that for some $\delta \in (0,1]$, $c_{2+\delta}(D)>0$ and any $D>0$
\begin{equation}\label{12}
{\sf E}|X_0(D)|^{2 + \delta} \leq c_{2+\delta}(D).
\end{equation}
Assume that condition \eqref{10} holds for all fields $X(D)$
with the same $\lambda > 4d(1+\delta)/\delta$ and $c_0$.
Then there exists $\nu = \nu(d,\lambda,\delta)>0$ such that for each $D>0$ and any $n \in \mathbb{N}$
\begin{equation}\label{13}
\sup_{x \in \mathbb{R}} |{\sf P}\bigl((S(U_n,D) - |U_n|{\sf E}X_0(D))/(\sigma(D) |U_n|^{1/2})
\leq x\bigr) - \Phi(x)|
\leq A|U_n|^{-\nu}
\end{equation}
where $A = A_0(d,\lambda)\max\{1, c_0{\sf E}X_0^2(D)\}\max\{1,
c_{2+\delta}(D)/\sigma^{2+\delta}(D)\}$ and \begin{equation}\label{14}
\sigma^2(D) = \sum_{j\in \mathbb{Z}^d} cov(X_0(D),X_j(D)) \ne 0.
\end{equation}
\end{theo}

{\bf Proof} is based on the classical blocks technique initiated by Bernstein, so we
only indicate the main steps and concentrate in the next Section on a
statistical version of this result.

For every $n \in \mathbb{N}$ introduce ${\sf
p}={\sf p}(n)= [n^{\alpha}]$ and ${\sf q}={\sf q}(n) = [n^{\beta}]$ where
$0<\beta <\alpha <1$ and $[\cdot]$ stands for an integer part of a number.
Consider $k =k(n)= [(2n+1)/(2{\sf p} + {\sf q})]$.  Then one can write
$[-n,n] =$ $I_1\cup I_1' \cup \ldots I_k\cup I_k' \cup I_k''$ where $I_m, I_m', I_k''$ are
disjoint intervals of the form $I_m = [a_m, a_m+2{\sf p}]$, $I_m' = (a_m+2{\sf
p},a_m+2{\sf p}+{\sf q})$ and $I_k'' = (a_k+2{\sf p}+{\sf q},n]$ ($I_k''$ can
be empty, $a_1 = -n, a_m \in [-n,n]$, $m=1,\ldots,k)$.

Set $B_i = I_{i_1} \times \ldots \times I_{i_d} \cap \mathbb{Z}^d$ where $i = (i_1,\ldots,i_d) \in
M_n = \{1,\ldots,k\}^d$ and let $V_n = \cup_{i \in M_n} B_i$.

It is easy to verify that for all n large enough
\begin{equation}\label{14ab}
\frac{1}{|U_n|}{\sf E}(S(U_n)- {\sf E} S(U_n) - (S(V_n)-{\sf
E}S(V_n)))^2 \leq \left(1 -
\frac{|V_n|}{|U_n|}\right)\sum_{j\in \mathbb{Z}^d}|cov(X_0,X_j)| \leq 4d
 n^{-\gamma}v(D)
\end{equation}
where $\gamma = \min\{1-\alpha,\alpha -\beta\}$ and the
series
\begin{equation}\label{14a}
\sum_{j\in \mathbb{Z}^d}|cov(X_0(D),X_j(D))| = v(D)
\end{equation}
converges in view of
\eqref{10} for $\lambda >d$.  Consequently,
\begin{equation}\label{15}
{\sf
E}\left|\frac{S(U_n) - |U_n|{\sf E}X_0}{\sigma(D)|U_n|^{1/2}} - \frac{S(V_n) -
|V_n|{\sf E}X_0}{\sigma(D)|U_n|^{1/2}}\right|\leq 2(d v(D))^{1/2}
n^{-\gamma/2}\sigma^{-1}(D).
\end{equation}

\begin{lem}\label{1} Let $X(D)=\{X_j(D), j \in \mathbb{Z}^d\}$ be a wide-sense
stationary random field such that \eqref{10} holds with some $\lambda > d$.
Then for all $n \in \mathbb{N}$ \begin{equation}\label{16} \left|\frac{var
S(U_n)}{|U_n|} - \sigma^2(D)\right| \leq a f(n,d,\lambda) \end{equation} where
$a = a_0(d,\lambda)c_0{\sf E}X_0^2(D)$ and $$ f(n,d,\lambda) =
\left\{\begin{array}{ll} n^{d-\lambda},&d<\lambda<d+1,\\ (1 +\ln n)/n,&\lambda
= d+1,\\ n^{-1},&\lambda > d+1.\end{array}\right.  $$ \end{lem} {\bf Proof.}
One has $$ |varS(U_n) - \sigma^2(D) |U_n|| \leq \sum_{i \in U_n}\sum_{j\notin
U_n} |cov(X_i,X_j)| = T_1 + T_2 $$ where the sums $T_1$ and $T_2$ are taken
respectively over the sets $\{i \in U_n, j \notin U_n,\|i-j\|\}=r, r >n\}$ and $\{i
\in U_n, j \notin U_n,\|i-j\|\}=r,\, r \in \{1,\ldots,n\}\}$.
Evidently,
$$
T_1
\leq 2d 3^{d-1}|U_n| c_0 {\sf E}X_0^2 \sum_{r>n} r^{d-1-\lambda} \leq 2d
 3^{d-1} c_0{\sf E}X_0^2n^{d-\lambda} |U_n| /(\lambda - d),
$$
$$
T_2 \leq
\sum_{r=1}^n\;\;\sum_{n-r <\|i\|\leq n, \|j-i\|=r} |cov(X_i,X_j)|
$$
$$
\leq
2dc_0 {\sf E}X_0^2 \sum_{r=1}^n r^{-\lambda}((2n+1)^d - (2(n-r)+1)^d)(2r+1)^{d-1}
\leq
4d^2 3^{d-1}c_0(2n+1)^{d-1}{\sf E}X_0^2 \sum_{r=1}^n r^{d-\lambda}.
$$
Using a trivial estimate
$$
\sum_{r=1}^n r^{-\mu} \leq 1 + \int_1^n x^{-\mu}dx, \;\;\mu >0,
$$
we come to relation \eqref{16}. The Lemma is proved.

Set $Y_j=S(B_j)- {\sf E}S(B_j)$ where
$B_j$ belongs to a collection of "large" blocks,
$j \in M_n=\{1,\ldots,k \}$, $k=k(n)$.
Clearly $Y_j = Y_j({\sf p}_n,D)$.
Introduce independent copies $Z_j$, $j \in M_n$,
of random variables $Y_j$, $j \in M_n$.

Then it is easily seen that for any $t \in \mathbb{R}$, $i^2=-1$ and all $n$
large enough $$ |{\sf E}\exp\{it\sum_{j \in M_n}Y_j\} - {\sf
E}\exp\{it\sum_{j\in M_n} Z_j\}| \leq 4|M_n| c_0 {\sf q}^{-\lambda} (2{\sf p}+1)^d|U_n|\leq 4
c_0 3^d n^{2d} {\sf q}^{-\lambda},
$$ $$ {\sf E}\left(
\frac{\sum_{j\in
M_n}Y_j}{\sigma(D)|U_n|^{1/2}}
\right)^2 \leq v(D)/\sigma^2(D),\;\;\; {\sf
E}\left(
\frac{\sum_{j\in
M_n}Z_j}{\sigma(D)|U_n|^{1/2}}
\right)^2 \leq
v(D)/\sigma^2(D).
$$
Thus
\begin{equation}\label{17}
\left|{\sf E}\exp\left\{
\frac{it\sum_{j\in M_n}Y_j}{\sigma(D)|U_n|^{1/2}}\right\}
 - {\sf E}\exp\left\{
\frac{\sum_{j\in M_n}Z_j}{\sigma(D)|U_n|^{1/2}}\right\}\right|
\leq \min\{4 c_0 3^d n^{2d} {\sf q}^{-\lambda}, 2|t|\sqrt{v(D)}/\sigma(D)\}.
\end{equation}

Using Lemma 1 one can verify that for all $n$ large enough
$$
{\sf E}\left(\frac{\sum_{j\in M_n} Z_j}{\sigma(D)|U_n|^{1/2}} - \frac{\sum_{j\in M_n} Z_j}
{(var\sum _{j\in M_n} Z_j)^{1/2}}\right)^2 \leq (4d n^{-\gamma} +
a\sigma^{-2}(D)f({\sf p},d,\lambda))^2
$$ where $\gamma$, $a$ and $f$ are the same as in
\eqref{14ab} and \eqref{16}.  Therefore for $\lambda > d+1$ we can write
\begin{equation}\label{18}
{\sf E}\left|\frac{\sum_{j\in M_n} Z_j}{\sigma(D)|U_n|^{1/2}} - \frac{\sum_{j\in M_n} Z_j}
{(var\sum _{j\in M_n} Z_j)^{1/2}}\right| \leq C_1 n^{-\tau}
\end{equation}
where $C_1 = C_1(d,a)\max\{1, \sigma^{-2}(D)\}$, $\tau = \min\{\gamma,
\alpha\}$.

Now the Esseen inequality implies that for every $T >0$ one has
\begin{equation}\label{19}
{\sf P}\left(\frac{S(U_n,D)}{\sigma(D)|U_n|^{1/2}} \leq x\right) - \Phi(x)|
\leq a_1 \int _{|t| \leq T} \left|\frac{{\sf
E}\exp\{it\frac{S(U_n)}{\sigma(D)|U_n|^{1/2}}\} -
\exp\{-\frac{t^2}{2}\}}{t}\right|dt
+a_2T^{-1}
\end{equation}
where $a_1$ and $a_2$ are absolute positive constants.

Applying
the Berry--Esseen estimate of the convergence rate in the CLT for independent summands
$Z_i, i \in M_n$, with finite absolute moments of order $2+\delta$ (see,
e.g., \cite{CT}, p. 322), using \eqref{15} -- \eqref{18} and estimating the integral in the
right hand side of \eqref{19} as a sum of integrals $\int_{|t|\leq 1/T}$ and
$\int_{1/T <|t|\leq T}$ and finally taking $T = bn^{\zeta}$ with appropriately small
$\zeta$ and specified $b>0$ we arrive at \eqref{13}. This completes the
proof of Theorem 1.

{\bf Remark 2.} There are many versions of the CLT for random fields
under various dependence conditions (see, e.g., \cite{BuZh},
\cite{Bolt}, \cite{BuMSU}, \cite{Dedec}, \cite{BulShash}) and
references therein). In the same manner we could use instead of the
maximal correlation coefficient $\rho$, e.g., the Rosenblatt-type
mixing coefficient. We proved here the CLT with rate because it
permits to establish the law of the iterated logarithm (announced in
\cite{BuLIL}) under the dependence conditions of the type
\eqref{10}. It is worth mentioning that to this end we need only
arbitrary slow power-type estimate of the convergence rate in the
CLT without specifying an exponent $\nu$ in \eqref{13}. We do not
provide here an explicit cumbersome expression for $\nu$. More
restrictive mixing conditions than \eqref{10}, i.e. $\rho(I,J) \leq
c_0 |I||J| \exp\{-a\; dist(I,J)\}$ where $a$ and $b$ are some
positive parameters, were recently used in \cite{SchSun},
\cite{SchSun1} (see also the references therein) for CLT and LIL. We
do not consider here growing subsets $U_n \subset \mathbb{Z}^d$ more
general than "integer" cubes. For generalizations of this kind we
refer to \cite{Bolt}, \cite{BulShash}.

\vskip1cm
{{\large \bf 4. Statistical version of the CLT}}
\vskip0,5cm

There are two ways for applications of Theorem 1. Namely, if we believe in the model describing the
stochastic behaviour of each FSU (see Section 1) then we can calculate ${\sf E}X_0$. However, the
problem for dependent FSUs is the following one. Now we cannot claim (in general) that the variance
of the sum $S(U_n,D)$ is equal to the sum of variances of summands. Thus for every $D>0$, in
contrast to the CLT for the Bernoulli scheme, i.e.
$$
\frac{S_n(D) - np(D)}{\sqrt{np(D)(1-p(D))}} \stackrel{{\mathcal D}}\to
Z\sim N(0,1)\;\;\;{\mbox as}\;\;n \to \infty
$$
discussed in Section 1 (here
$p(D) = {\sf E}X_0$, $Z$ is a standard normal r.v.), the relation
\begin{equation}\label{19a}
\frac{S(U_n,D) - |U_n|{\sf E}X_0(D)}{\sigma(D)|U_n|^{1/2}}
\stackrel{{\mathcal D}}\to Z\sim N(0,1)\;\;\;{\mbox as}\;\;n \to \infty
\end{equation}
contains an unknown function $\sigma(D)$. As usual
$\stackrel{{\mathcal D}}\to$ stands for weak convergence
of random variables distributions.

Fortunately it is possible to construct a sequence of nonnegative
statistical estimates $\widehat{C}(U_n,D)$
for $\sigma^2(D)$ such that for any $D>0$
\begin{equation}\label{20}
\widehat{C}(U_n,D) \stackrel{{\sf P}}\to \sigma^2(D)\;\;\;{\mbox as} \;\;n\to \infty
\end{equation}
where $\stackrel{{\sf P}}\to$ means the convergence in probability as usual. We employ here a
family of consistent statistical estimates introduced in \cite{BulVr} for random fields \footnote
{Vector-valued random fields satisfying other dependence conditions are studied in \cite{Bul}.},
for stochastic processes we refer to the paper \cite{PelSh}.

Then by virtue of \eqref{19a} and \eqref{20} we come, for every
$D>0$ (if $\sigma(D)^2 \ne 0$), to the formula
\begin{equation}\label{21} (\widehat{C}(U_n,D)|U_n|)^{-1/2}(S(U_n,D)
- |U_n|{\sf E}X_0(D)) \stackrel{{\mathcal D}}\to Z\sim
N(0,1)\;\;\;{\mbox as}\;\;n \to \infty.
\end{equation}

In other words a random normalization is used in the CLT.

Consequently to determine (approximately) for a given value $\gamma \in (0,1)$
the threshold $x_{\gamma}$ we can apply the following analogue of formula
\eqref{3} \begin{equation}\label{22} {\sf P}(S_n \geq x) \approx 1 - \Phi(\xi)
\end{equation}
where $\xi = (x-n{\sf E}X_0(D))/\sqrt{\widehat{C}(U_n,D)|U_n|}$, $x \in \mathbb{R}$, $n \in
\mathbb{N}$.  However, now in the right hand side of \eqref{22} there is a r.v. $\Phi(\xi)$, i.e.
we use $\Phi(\xi)$ as statisitical estimate for ${\sf P}(S_n <x)$. Note that we have used only the
value ${\sf E}X_0(D)$ provided by the model of stochastic behaviour for FSUs and we did not suppose
here that the collective effect of the evolution of cells under irradiation is described by
independent binary random variables.

Another way of using Theorem 1 is to construct approximate confidence
intervals for the unknown mean value ${\sf E}X_0(D)$ without hypotheses
concerning the explicit formulas (discussed in Section 1) for distribution of
random variables $X_j$, $j \in U_n$.

Thus in both cases it is desirable to establish the CLT for dependent random fields
using random normalization.

{\bf Remark 3.}
As far as we know, in previous applications of the CLT to NTCP models for
independent FSUs the question of convergence rate was not raised, so Section
1 covers this gap. However, the same question in case of dependent FSUs is
more involved. We intend to investigate
the accuracy of the proposed
model in a special publication. One can consider Theorem 1 as the first step
in this direction. Moreover, we can obtain the power-type estimate in the CLT with random
normalization. However, the rate of convergence will be slower than that
for independent random summands.  The effect of convergence rate in the CLT
sensitivity to the dependence conditions was demonstrated for positively or
negatively associated random fields in \cite{BuASS}.

For $j \in U \subset \mathbb{Z}^d$ $(1 \leq |U| < \infty)$ and $b =
b(U) > 0$  set
\begin{equation}\label{23}
K_j(b) = \{t \in \mathbb{Z}^d: \|j - t\| \leq b\},\;\; Q_j =
Q_j(U,b) = U \cap K_j(b),
\end{equation}
\begin{equation}\label{24}
\widehat{C}(U,D) = \frac{1}{|U|}\sum_{j \in U} |Q_j|
\left(\frac{S(Q_j,D)}{|Q_j|} -
\frac{S(U,D)}{|U|}\right)^2.
\end{equation}

Note that
the averaged variables $S(Q_j,D)/|Q_j|$ arise
for dependent summands (in contrast to the traditional
estimates of variance used for independent observations).

\begin{theo}\label{2}
Let the conditions of Theorem \ref{1} be satisfied.
Let $\{U_n\}_{n \in \mathbb{N}}$ be a sequence of "integer" cubes, i.e.
$U_n = [-n,n]^d \cap \mathbb{Z}^d$, $n\in \mathbb{N}$ $(d\geq 1)$.
Assume $b(U_n) = b_n$ where
$\{b_n\}_{n \in \mathbb{N}}$ is a sequence of positive integers
such that
\begin{equation}\label{25}
b_n \to \infty, \;\; b_n = o(n)\;\;{\mbox as}\;\;n \to \infty.
\end{equation}
Then for every $D>0$ relation \eqref{21} holds if $\sigma(D) \ne 0$.
\end{theo}

{\bf Proof}. The estimate $\widehat{C}(U,D)$ introduced by means of
\eqref{24} and  $\sigma^2(D)$ are invariant under the transformation
$X_j \mapsto X_j - {\sf E}X_0$, $j \in \mathbb{Z}^d$. So, without
loss of generality we can further on assume that ${\sf E}X_0 = 0$.
Let $\|\xi\|_L$ stand for the norm of a real-valued random variable
$\xi$ in a space $L = L^1(\Omega,\mathcal{F},{\sf P})$. For any $U_n
\subset \mathbb{Z}^d$ and $n \in \mathbb{N}$ one has
$$
\|\widehat{C}(U_n,D) - \sigma^2(D)\|_L \leq I_1(U_n,D) + I_2(U_n,D) +
I_3(U_n,D)
$$
where
$$
I_1(U_n,D) = \frac{1}{|U_n|}\left\|\sum_{j \in
U_n}|Q_j|\left\{\left(\frac{S(Q_j)}{|Q_j|} -
\frac{S(U_n)}{|U_n|}\right)^2
 -
\left(\frac{S(Q_j)}{|Q_j|}\right)^2\right\}\right\|_L,
$$
$$
I_2(U_n,D) = \frac{1}{|U_n|}\left\|\sum_{j \in U_n}
\frac{1}{|Q_j|}\Bigl(S^2(Q_j) -
{\sf E}S^2(Q_j)\Bigr)\right\|_L,
$$
$$
I_3(U_n,D) = \left|\frac{1}{|U_n|}\sum_{j\in U_n}
\frac{1}{|Q_j|}{\sf E}S^2(Q_j) - \sigma^2(D)\right|.
$$
Here and below $S(U) = S(U,D)$ for $U \subset \mathbb{Z}^d$ and $D>0$. We have
\begin{equation}\label{27a}
|Q_j|^{-1} ES^2(Q_j) \leq
v(D),\;\;j \in \mathbb{Z}^d,
\end{equation}
where
$v(D)$ is the same as in \eqref{14a}.
By virtue of condition \eqref{25} it is clear that
$$
I_1(U_n,D) \leq |U_n|^{-3}{\sf E}S^2(U_n)\sum_{j \in U_n} |Q_j|
+\; 2|U_n|^{-2}\sum_{j\in U_n}{\sf E}|S(Q_j)S(U_n)|
$$
\begin{equation}\label{26}
\leq v(D)\{|K_0(b_n)||U_n|^{-1} +
2|K_0(b_n)|^{1/2}|U_n|^{-1/2}\} \to 0\;\;\mbox{as}\;\;n \to \infty.
\end{equation}
For a fixed $c > 0$ introduce the functions
\begin{equation}\label{27}
h_1(x) = \mbox{{\rm sign}}(x)\min\{|x|,c\},\;\; h_2(x) = x -
h_1(x), \;\;x \in \mathbb{R}.
\end{equation}
Given a nonempty finite set $Q \subset \mathbb{Z}^d$ let
$$\overline{S}(Q) = S(Q)/\sqrt{|Q|}.$$

Note that
\begin{equation}\label{28}
I_2(U_n,D) \leq \sum_{p,m=1}^2 I_2^{(p,m)}(U_n,D)
\end{equation}
where
$$
I_2^{(p,m)}(U_n,D) = \frac{1}{|U_n|}\Bigl\|\sum_{j \in U_n}
h_p(\overline{S}(Q_j))h_m(\overline{S}(Q_j)) -
{\sf E}h_p(\overline{S}(Q_j))
h_m (\overline{S}(Q_j))\Bigr\|_L.
$$
For $b,n \in \mathbb{N}$ introduce the sets
$$
T_n^{(b)} = \{s \in U_n: \inf_{t \in \partial U_n} \|s - t\| \leq
b\}
$$
where $\partial U_n = \{j \in U_n: \exists q \notin U_n\;\;\mbox{such
that}\;\; \|j-q\|=1\}$.  Put $T_n = T_n^{(b_n)}$, $n \in \mathbb{N}$, where
$b_n$ meet condition \eqref{25}.  Due to \eqref{27a} one has $$ I_2^{(1,2)}
(U_n,D) + I_2^{(2,1)}(U_n,D) + I_2^{(2,2)}(U_n,D) $$ $$ \leq 2|U_n|^{-1}
\sum_{j \in U_n} \biggl( 2{\sf E}| h_1(\overline{S}(Q_j))h_2(\overline{S}(Q_j))
| + {\sf E} h^2_2(\overline{S}(Q_j)) \biggr) $$ $$ \leq 2 \biggl( 2{\sf E}|
h_1(\overline{S}(K_0(b_n)))h_2(\overline{S}(K_0(b_n)))|
+ {\sf E}
h_2^2(\overline{S}(K_0(b_n)))
+
\;3|T_n||U_n|^{-1}v(D)\biggr)
$$
$$
\leq
4\Bigl(v(D){\sf E} \bigl( \overline{S}^2(K_0(b_n))
{\bf 1}
\bigl\{|\overline{S}(K_0(b_n))|
\geq
c\bigr\}\bigr)\Bigr)^{1/2}
$$
$$
+\;
2 {\sf E}\bigl(\overline{S}^2(K_0(b_n)
{\bf 1}
\bigl\{|\overline{S}(K_0(b_n))|
\geq
c\bigr\}\bigr)
$$
\begin{equation}\label{28a}
+\; 6|T_n||U_n|^{-1}v(D)
\end{equation}
where ${\bf 1}$ is an indicator function.

Condition \eqref{27a} implies that
\begin{equation}\label{28b}
|T_n||U_n|^{-1} \to 0\;\;\mbox{as}\;\;n\to \infty.
\end{equation}

It is easy to see that a family
$\{\overline{S}^2(K_0(b_n))\}_{n=1}^{\infty}$ is uniformly
integrable. Consequently, taking into account \eqref{28a}, for any
$\varepsilon > 0$ we can find $c=c(\varepsilon)$ such that for all $n$ large
enough \begin{equation}\label{29} I_2^{(1,2)} (U_n,D) + I_2^{(2,1)}(U_n,D) +
I_2^{(2,2)}(U_n,D) < \varepsilon,  \end{equation} furthermore,
\begin{equation}\label{30}
(I_2^{(1,1)}(U_n,D))^2 \leq
|U_n|^{-2}\sum_{j,t \in U_n} \Bigl| \,\mbox{{\rm cov}}
\biggl(
h_1^2(\overline{S}(Q_j)),
h_1^2(\overline{S}(Q_t))
\biggr)\Bigr|.
\end{equation}
In view of \eqref{27} we obtain the inequalities
$$
|U_n|^{-2}\sum_{j,t \in U_n, \|j-t\|\leq 4b_n} \left|\,\mbox{{\rm
cov}}\biggl(
h_1^2(\overline{S}(Q_j)),
h_1^2(\overline{S}(Q_t))\biggr)
\right|
$$
$$
\leq 2c^2|U_n|^{-2} \sum_{j,t\in U_n,\|j-t\|\leq 4b_n}
{\sf E}h_1^2(\overline{S}(Q_j))
$$
\begin{equation}\label{32}
\leq 2^{2d+1}c^2|U_n|^{-1}|K_0(b_n)|v(D).
\end{equation}
Now condition \eqref{10} with $\lambda > d$ entails the estimate
$$
|U_n|^{-2}\sum_{j,t \in U_n, \|j-t\| > 4b_n} \left|\,\mbox{{\rm
cov}}\left(
h_1^2(\overline{S}(Q_j)),
h_1^2(\overline{S}(Q_t))\right)
\right|
\leq c_0 c^4|U_n|^{-2}\sum_{j,t \in U_n, \|j-t\| > 4b_n}
|Q_j||Q_t| \|j-t-2b_n\|^{-\lambda}
$$
\begin{equation}\label{33}
\leq c_0c^4 d4^d |U_n|^{-1}|K_0(b_n)|^2 \sum_{r>2b_n} r^{d-1-\lambda}.
\end{equation}
Taking into account \eqref{29} -- \eqref{33}, \eqref{25} and \eqref{10} with
$\lambda \geq 2d$, we verify that \begin{equation}\label{34} I_2(U_n,D) \to
0\;\;\mbox{as}\;\; n\to \infty.  \end{equation} Now observe that
$$ |U_n|^{-1}
\sum_{j\in U_n} |Q_j|^{-1} {\sf E}S^2(Q_j) = |U_n|^{-1}|U_n\setminus
T_n||K_0(b_n)|^{-1} {\sf E}S^2(K_0(b_n)) $$ $$ + |U_n|^{-1}\sum_{j\in T_n}
|Q_j|^{-1} {\sf E}S^2(Q_j).  $$ According to \eqref{16} the following relation is
valid $$ |K_0(b_n)|^{-1}{\sf E}S^2(K_0(b_n)) \to \sigma^2(D)\;\;\mbox{as}\;\; n
\to \infty.  $$
Due to  \eqref{25} and \eqref{27a} we conclude that
\begin{equation}\label{35}
I_3(U_n,D) \to 0\;\;\mbox{as}\;\;n\to \infty.
\end{equation}
Relations \eqref{26}, \eqref{34} and \eqref{35} yield \eqref{19}.
The proof of Theorem 2 is complete.

\vskip1cm
{\large {\bf 5. Concluding Remarks}}

\vskip0,3cm

Besides the concept of a functional reserve of an organ undergoing irradiation it seems desirable
to use the models taking into account the geometrical configuration of survived FSUs (or/and
cells). Here the concepts of random clusters appear naturally. In this regard we refer to a quite
recent paper \cite{Thames} (cf. \cite{Klepp}). The stochastic models of disordered structures
(involving point random fields) could be applied also to describe the damage volumes of irradiated
organ. Note that it is possible to describe the relations between various FSUs (cells) in terms of
random graphs and study the models where some vertices (or edges) are destroyed at random.  An
interesting problem is to find the optimal dose of irradiation taking into account not only the
complication probabilities but the balance of conditions for irradiated organ (tissue) and its
normal environment.  Moreover, it is important to consider non uniform irradiation, another problem
is to study a population of non-identical patients (see, e.g., \cite{Yorke01},\cite{FN}).  To
conclude we mention a deep problem of constructing dynamical models describing the evolution of an
irradiated organ in space and time.

\vskip1cm
{\large{\bf  Acknowledgements}}

\vskip0,3cm The authors thank Professor L.Klepper for drawing their attention to stochastic models
in Radiobiology. A.Bulinski is grateful for invitation to the V\"axj\"o University where in August
2004 the results of this joint paper were obtained and he is grateful also to the Dept. of
Probability Theory and Stochastic Models of the University Paris-VI where the final version of the
paper was written.


\vskip1cm

\begin{thebibliography}{99}

\bibitem{Bolt} Bolthausen E., {\it On the central limit theorem for stationary mixing random fields},
Ann. Probab., {\bf 10}(1982), 1047--1050.

\bibitem{BuZh} Bulinski A.V., Zhurbenko I.G., {\it A central limit theorem for additive random functions},
Theory Probab. Appl., {\bf 21}(1976), 687--697.

\bibitem{BuMSU} Bulinski A.V., {\it Limit Theorems under Weak Dependence Conditiions}, MSU, 1989
(in Russian).

\bibitem{BuASS} Bulinski A.V., {\it On the convergence rates in the CLT for positevely or negatively
dependent random fileds}, In.: I.A.Ibragimov and A.Yu.Zaitsev (Eds.) Proc. of the Kolmogorov
Semester Int. Euler Math. Inst., St.-Petersburg (02.03.1993 -- 23.04.1993). Gordon and Breach,
1996, 3 -- 14.

\bibitem{BulVr} Bulinski A.V., Vronski M.A., {\it Statistical variant of the
central limit theorem for associated random fields}, Fund. Appl. Math., {\bf 2}(1996), 891--908 (in
Russian).

\bibitem{Bul} Bulinski A.V., {\it Statistical version of the central limit
theorem for vector-valued random fields}, Math. Notes, {\bf
76}(2004), 455--464.

\bibitem{BulShash} Bulinski A.V., Shashkin A.P., {\it Rates in the central limit hteorem for weakly
dependent random variables}, J. Math. Sci., {\bf 122}(2004), 3343--3358.

\bibitem{BuLIL} Bulinski A.V., {\it A law of the iterated logarithm for random fields with power decay of
correlations}, Rev. Appl. Indust. Math., {\bf 11}(2004), 503--504 (in Russian).

\bibitem{B} Bulinski A.V., {\it Stochastic models in micro- and macroworld},
Proc. Semin. "Time, Chaos and Mathematical Problems",  MSU, {\bf 3}(2004), 231--254 (in Russian).

\bibitem{CT} Chow S.C., Teicher H., {\it Probability Theory. Independence, Interchangebility,
Martingales}, Springer, New York etc., 1997 (3-d ed).

\bibitem{Dedec} Dedecker J., {\it A central limit theorem for stationary random
fields},  Probab. Th. Rel. Fields, {\bf 110}(1998), 397--426.

\bibitem{Dob} Dobrushin R.L., {\it Description of a random field by means of
conditional probabilities and conditions of its regularity}, Theory
Probab. Appl., {\bf 13}(1968), 197--225.

\bibitem{Douk} Doukhan P., {\it Mixing: Properties and Examples}, Lecture Notes in Statistics,
{\bf 85}(1994), Springer--Verlag.

\bibitem{FN} Fenwick J.D., Nahum A.E., {\it Series model volume effects in a
population of non-identical patients: how low is low?}, Phys. Med. Biol., {\bf 46}(2001),
1815--1834.

\bibitem{Jack} Jackson A., Kutcher G.J., Yorke E.D., {\it Probability of
radiation-induced complications for normal tissues with parallel architecture subject to
non-uniform irradiation}, Med. Phys., {\bf 20}(1993), 613--625.

\bibitem{JonDale} Jones B., Dale R.G., {\it Mathematical models of tumours and normal tissue response},
Acta Oncol., {\bf 38}(1999), 883--893.

\bibitem{Klepp} Klepper L.I., {\it The occurrence of radiation complications in
body organs and tissues}, Med Tekh., {\bf 5}(2000), 36--40 (in Russian).

\bibitem{LMHK} Lind B.K., Mavroidis P., Hyodynmaa S., Kappas C., {\it Optimization
of the dose level for a given treatment plan to maximize the complication-free tumor cure}, Acta
Oncol., {\bf 38}(1999), 787--798.

\bibitem{MBD} Moiseenko V., Battista J., Van Dyk J.,
{\it Normal tissue complication probabilities: dependence on choice of biological model and
dose-volume histogram reduction scheme}, Int. J. Radiat. Oncol. Biol. Phys., {\bf 46}(2000),
983--993.

\bibitem{Nie} Niemierko A., Goitein M., {\it Modeling of normal tissue response
to radiation: the critical volume model}, Int. J. Radiat. Oncol. Biol. Phys., {\bf 25}(1993),
135--145.

\bibitem{PelSh} Peligrad M., Shao Q-M., {\it Estimation of variance of partial sums of an associated
sequence of random variables}, Stoch. Proc. Appl., {\bf 56}(1995), 307--319.

\bibitem{Pet} Petrov V.V., {\it Limit Theorems
of Probability Theory: Sequences of Independent Random Variables}, Oxford Univ. Press, 1995.

\bibitem{PPVSK} Philippsens M.E., Pop L.A., Visser A.G., Schelleksns S.A.,
Van Der Kogel A.J., {\it  Dose-volume effects in rat thoracolumbar spinal chord: an evaluation of
NTCP models}, Int. J. Radiat. Oncol. Biol. Phys., {\bf 60}(2004), 578--590.

\bibitem{SchSun} Schmuland B., Sun W., {\it A central limit theorem and law of the iterated logarithm for
random fields with exponential decay of correlations}, Canad. J. Math. {\bf 56}(2004), 209--222.

\bibitem{SchSun1} Schmuland B., Sun W.A., {\it The law of large numbers and law or the iterated logarithm
for infinite dimensional intaracting diffusion processes}, Inf. Dim.
Anal., Quantum Probab. Rel. Topics (IDAQP), {\bf 6}, 489--503.

\bibitem{Smith} Smith A.R. (Ed.), {\it Medical Radiology: Radiation Therapy Physics},
Springer -- Verlag, Berlin, 1995.

\bibitem{SNSG} Stavrev P., Stavreva N., Niemierko A., Goitein M.,
{\it Generalization of a model of tissue response to radiation based on the idea of functional
subunits and binomial statistics}, Phys. Med. Biol., {\bf 46}(2001), 1501--1518.

\bibitem{SSWF} Stavreva N., Stavrev P., Warkentin B., Fallone B.G., {\it Derivation
of the expressions for gamma 50 and D50 for different individual TCP and NTCP models}, Phys. Med.
Biol., {\bf 47}(2002), 3591--3604.

\bibitem{Thames} Thames H.D., Zhang M., Tucker S.L., Liu H.H., Dong L., Mohan
R., {\it Cluster models of dose-volume effects}, Int. J. Radiat. Oncol. Biol. Phys., {\bf
59}(2004), 1491--1504.

\bibitem{LDSS} van Luijk P., Delvigne T.C., Schilstra C., Schippers J.M., {\it Estimation of parameters
of dose-volume models and their confidence limits}, Phys. Med. Biol., {\bf 48}(2003), 1863--1884.

\bibitem{Weiss} Weisstein E.W., {\it "de Moivre--Laplace Theorem"}.
From {\it MathWorld} -- A Wolfram Web
Resource. http://mathworld.wolfram.com/deMoivre-LaplaceTheorem.html

\bibitem{Yorke01} Yorke E.D., {\it Modeling the effect of inhomogeneous dose distributions in normal
tissues}, Semin. Radiat. Oncol., {\bf 11}(2001), 197--209.



\end{thebibliography}
\end{document}